\documentclass[a4paper]{jpconf}

\usepackage{graphicx}

\usepackage{amsmath}
\usepackage{amsthm}
\usepackage{amssymb}
\usepackage{paralist}
\usepackage{graphics}
\usepackage{epsfig}
\usepackage[colorlinks=true]{hyperref}
\hypersetup{urlcolor=blue, citecolor=red}
\usepackage{algorithmic}
\usepackage{url}
\usepackage[section]{placeins}

\newtheorem{definition}{Definition}

\begin{document}
\title{Mathematical Modeling for Tomography in Domains with Reflecting Obstacles}

\author{Kamen Lozev}

\address{3161 S. Sepulveda Boulevard, Apt. 307, Los Angeles, CA 90034}

\ead{kamen@ucla.edu}

\begin{abstract}
This work develops new numerical methods for the solution of the tomography problem in domains with reflecting obstacles. We compare the solution's performance 
for Lambertian reflection, for classical tomography with unbroken rays and for specular reflection. Our numerical method using Lambertian reflection improves 
the solution's accuracy by an order of magnitude compared to classical tomography with unbroken rays and for tomography in the presence of a specularly reflecting 
obstacle the numerical method improves the solution's accuracy approximately by a factor of three times. We present efficient new algorithms for the solution's software 
implementation and analyze the solution's performance and effectiveness. The new method from this work for reducing the number of equations in tomography linear systems
is applicable to improving the performance of a wide class of algebraic tomographic reconstruction methods.   
\end{abstract}

\section{Introduction}

\bigskip

Let $f(x)$ be a continuous function in $\Omega$, where $\Omega=\Omega_0 \backslash  \overline{\Omega_1}$ and $\Omega_0$ is a 
compact convex set in $\mathbb{R}^2$ with a smooth boundary and $\Omega_1$ a convex obstacle with a smooth boundary such that
$\overline{\Omega_1} \subset \Omega_0 \subset \mathbb{R}^2$. 
Consider $\partial{\Omega_0}$ as the observation boundary of $\Omega_0$ and the points of this boundary 
as transmitters and receivers of ray signals or ray solutions of the wave equation 
\begin{equation}\label{waveequation}
u_{tt} - c^2(x)\Delta{u}=0
\end{equation}

where $c(x)>0$ is the variable speed of sound in $\Omega$  and $ u |_{\partial{\Omega_1}} = 0 $.     

The data for the tomography problem are all integrals $\int_{\gamma} f(l) dl = C_{\gamma}$ where $\gamma$ are  
rays in $\Omega$ that have both of their end points in $\partial{\Omega_0}$. One end point is a transmitter 
and the other end point of each of the rays a receiver. The classical tomography problem is to find $f(x)$ in 
$\Omega$ knowing the values of all integrals $C_{\gamma}$ where $\gamma$
are all straight line segments or unbroken rays in $\Omega$. 
In the case of a domain $\Omega_0$ without an obstacle $\Omega_1 \subset \Omega_0$, this problem is 
widely studied theoretically and numerically \cite{Na1,KS,NW,F}. When there are obstacles present, 
this problem is much less studied. Key theoretical work is done in \cite{E1,E2,E3,E4} for domains 
with one obstacle and with both broken rays, i.e. rays reflecting at the obstacle, and unbroken rays. 
A key result from \cite{E1} is that the tomography problem in the presence of one reflecting obstacle 
is well-posed. In other words, if an obstacle $\Omega_1 \subset \Omega_0$ is present then we have the well-posed problem of recovering 
$f(x)$  in $\Omega_0 \backslash \Omega_1$ from the set $C_{\gamma}$ where $\gamma$ are all broken and unbroken 
rays in the domain starting and ending at the observation boundary. This problem is called the Broken Ray Tomography Problem\cite{E1}. 

In a basic tomography setup transmitters and receivers of wave signals are placed at the domain's boundary $\partial{\Omega_0}$. 
Ray signals are generated by the transmitters and received by the receivers. Travel times for signal propagation from transmitters to 
receivers are measured and these travel times $T(A, B)$ are the values of line integrals of a function $f(x)=\frac{1}{c(x)}$ 
where $c(x)>0$ is the speed of sound at point $x \in \Omega_0 \backslash \Omega_1$. This measurement procedure gives the $C_{\gamma}$ 
data for solving the Tomography Problem by relating signal travel times to the values of line integrals of $f(x)$. 
Given sufficient data, $f(x)$ and from here the velocity $c(x)$ are computed with tomographic reconstruction algorithms \cite{Na1,KS,NW,F}. 
The tomographic algorithms presented in the next section require and compute the ray path $\gamma$ of each ray and, in order to compute these ray paths, 
we consider reflection models at $\partial{\Omega_1}$ and models of the speed of sound in $\Omega$. 

The theory of broken ray tomography from \cite{E1,E2,E3,E4} is based on a reflection model at $\partial{\Omega_1}$ that is mirror-like i.e. 
the angle of incidence is equal to the angle of reflection. In this paper, we consider a Lambertian reflection model in which incident rays are
reflected at the obstacle in all possible directions and present results that show that the broken ray tomography reconstruction 
error is smaller when we consider Lambertian reflection. 

We consider a mathematical model of the speed of sound 
\begin{equation}\label{close_to_constant_speed_of_sound}
c(x)=c_{o}+\epsilon(x)
\end{equation} 
where $x \in \Omega$ that models the speed of sound as a continuous function close to a known constant $c_{o}$. 
It is shown in \cite{RO} that for sufficiently small $\epsilon(x)$ waves propagate along the known geodesics of $c_{o}(x)$ 
when 
\begin{equation}\label{continuous_speed_of_sound}
c(x)=c_{o}(x)+\epsilon(x) 
\end{equation}
where $c(x)$, $c_{o}(x)$ and $\epsilon(x)$ are continuous functions in $\Omega$. The acoustic geodesics for constant speed of sound 
$c_o(x)=c_{o}$ are straight lines when there is no obstacle. Therefore, for the model \ref{close_to_constant_speed_of_sound}, we consider two cases. 
In the first case, $\gamma = \tilde{\gamma_1}$ is an unbroken ray composed of a straight line segment 
$\tilde{\gamma_1}$. In the second case of a broken ray, $\gamma=\tilde{\gamma_1} \bigcup \tilde{\gamma_2}$ is the union of two straight 
line segments that intersect at a reflection point at the obstacle. 

For unbroken rays the travel time or time of flight is
$$ T(A,B)=\int_{\gamma}\frac{ds}{c_o + \epsilon{(x(s))}} = \int_{\tilde{\gamma_1}}\frac{ds}{c_o + \epsilon{(x(s))}} $$ and this 
model leads to the classical tomography problem with $f(x)=\frac{1}{c_o+\epsilon(x)}$. 

\bigskip

For broken rays and a known obstacle, we know that the acoustic wave $u(x)$ propagates along the known straight line segments 
$\tilde{\gamma_1}$ and $\tilde{\gamma_2}$. $\gamma_1$ and $\gamma_2$ are known because in addition to the time of flight, our data
measurement procedure gives the end points of the ray $\gamma$ and its initial velocity, which in turn imply the reflection point of 
$\gamma$ at the known obstacle $\Omega_1$. 
Then $$ T(A,B)=\int_{\gamma}\frac{ds}{c_o + \epsilon{(x(s))}}=\int_{\tilde{\gamma_1}}\frac{ds}{c_o + \epsilon{(x(s))}} + \int_{\tilde{\gamma_2}}\frac{ds}{c_o+\epsilon{(x(s))}} $$
and, when data of this type is added to the set of measurements for the time of flight for unbroken rays, this gives the set 
$C_{\gamma}$ for the Broken Ray Tomography Problem with $f(x)=\frac{1}{c_o+\epsilon(x)}$. 
In other words, the data set $C_{\gamma}$ for the Broken Ray Tomography Problem contains the travel times of all broken and unbroken 
rays in the domain that start and end at the observation boundary.
\bigskip

\section{Numerical Solution of the Broken Ray Tomography Problem for Lambertian Reflection}\label{section_BRTL}

The first algorithm for the Broken Ray Tomography Problem is presented in \cite{L1} 
for a known obstacle $\Omega_1$, specular reflection and the model \ref{close_to_constant_speed_of_sound} for the speed of sound. 
This work extends the first numerical solution of the Broken Ray Tomography Problem and develops an 
algorithm for finding the velocity structure of $\Omega$ for Lambertian reflectance at $\partial{\Omega_1}$. The broken ray tomography 
problem for an obstacle with Lambertian reflectance is solved with the following algorithm that constructs the finite set of broken and 
unbroken rays  and computes the associated ray travel times by numerical integration or from ray travel time data.

\begin{algorithmic}\label{algorithm_BRTL}
\REQUIRE Domain $\Omega_0$
\REQUIRE Obstacle $\Omega_1 \subset \Omega_0$
\REQUIRE Finite set of receiver points $R=(x_r, y_r)$ on $\partial{\Omega_0}$
\REQUIRE Finite set of transmitter points $T=(x_t, y_t)$  on $\partial{\Omega_0}$
\REQUIRE Finite set of obstacle boundary points $H=(x_h, y_h)$ on $\partial{\Omega_1}$
\REQUIRE Number of broken rays $n_b$. 
\REQUIRE Number of unbroken rays $n_u$.

\COMMENT{Algorithm for reconstructing $f(x,y)$ in $\Omega$ in the presence of obstacle $\Omega_0$}

\STATE{initialize an empty list L that will contain all broken and unbroken rays}

\FOR{$p = 1 \to n_b$}
  \STATE{generate a random broken ray r from input data R, T, H, $\Omega_0$ and $\Omega_1$}
  \COMMENT{a broken ray is a unique triple $((x_r, y_r), (x_t, y_t), (x_h, y_h))$ that does not intersect the obstacle except at the reflection point $(x_h, y_h)$}
  \STATE{add r to L}
\ENDFOR

\FOR{$p = 1 \to n_u$}
  \STATE{generate a random unbroken ray r from input data R, T, H, $\Omega_0$ and $\Omega_1$}
  \COMMENT{an unbroken ray is a unique pair $((x_r, y_r), (x_t, y_t))$ that does not intersect the obstacle}
  \STATE{add r to L}
\ENDFOR

\STATE{randomize L as a preprocessing step before starting the Kaczmarz method}

\FORALL{rays r in L}
\STATE{Compute travel time $p_r$ for ray r. In numerical simulations, this is obtained via nimerical integration in $\Omega_0$  along r. Store $p_r$ in list P.}  
\ENDFOR

\STATE{Reconstruct $f(x,y)$ by the Kaczmarz method \cite{K} for the linear system $Wf=P$ associated with rays r in L and corresponding travel times $p_r$ in P.}  

\end{algorithmic}

There are many algebraic reconstruction methods for obtaining a linear system of the form $Wf=P$ from a ray set $L$ and corresponding travel times $P$ 
where $W$ is a two dimensional matrix, f an unknown vector of the values of the function in each cell of a computation grid, and P the vector of 
travel times \cite{KS}. This linear system is solved in the above algorithm by the Kaczmarz method. Section \ref{drburzo} presents a new method for constructing
$Wf=P$ that reduces the number of equations of the linear system in order to reduce the run time of the algorithm.  
As input to the above algorithm we give $n_b$ and $n_u$ to be equal or approximately equal to the maximum number of broken and unbroken rays for the finite 
input sets R, T and H, or alternatively, during ray generation, generate new rays until all rays with end points in the input sets are generated. 
In order to approximate the requirements of the theory of broken ray tomography for inclusion of all broken and unbroken rays, we provide as input to the algorithm 
as many transmitter and receiver points on the observation boundary as possible so that each set approximates the set of all points on the observation boundary
and in order for the algorithm to approximate the set of all broken and unbroken rays.
The above algorithm uses Lambertian reflection by considering broken rays as random triples $((x_r, y_r), (x_t, y_t), (x_h, y_h))$ 
and uses all rays in the domain. 

\section{Experimental Results}

In order to show the effectiveness of the numerical solution of the broken ray tomography problem with Lambertian reflection, a Java implementation of the above
algorithm compares the reconstructed values of $f(x,y)$ in $\Omega_0$ with the known values of $f(x,y)$ for the same test function $f(x,y)=K\sqrt{(x-x_0)^2+(y-y_0)^2}$ 
, where $(x_0,y_0)$ is the center of the computation grid, and test environment with a square obstacle as in \cite{L1}. 
Table \ref{BRTLambertianData1} compares the reconstruction error for classical tomography without reflection and broken ray tomography with Lambertian reflection. 

\begin{table}[h]
\begin{tabular}{|c|cc|cc|}
	\hline
Experiment & ART Error  & ART Iterations & BRTL Error  & BRTL Iterations  \\
  \hline

1 & 1.269592e-004 & 54293 & 1.272610e-005 & 32635 \\ 
2 & 1.300814e-004 & 54204 & 5.272083e-006 & 40557 \\ 
3 & 1.478464e-004 & 45923 & 7.530267e-006 & 37251 \\ 
4 & 1.642026e-004 & 42480 & 1.110909e-005 & 33104 \\ 
5 & 1.927561e-004 & 23037 & 1.163345e-005 & 33424 \\ 
6 & 1.985439e-004 & 33705 & 1.896151e-005 & 28522 \\ 
7 & 8.991362e-005 & 88190 & 1.401511e-005 & 30061 \\ 
8 & 1.641201e-004 & 40717 & 1.955676e-005 & 29687 \\ 
9 & 1.089771e-004 & 64932 & 7.102756e-006 & 40126 \\ 
10 & 1.934379e-004 & 32899 & 2.149935e-005  & 28723 \\ 

  \hline

Average & 1.516838e-004 & 48038 & 1.294065e-005 & 33409 \\

  \hline
\end{tabular}
\caption{Error and number of iterations for broken ray tomography with Lambertian reflection(BRTL) at the boundary of the reflecting obstacle 
for a fixed number of 126050 rays. The average error for BRTL is 1.294065e-005 and the average number 
of iterations for finding a solution is 33409, and are shown in the right two columns of the table. 
The results for classical tomography with the ART method with the same number of rays are shown in the left two columns of the table. 
The average error for ART is 1.516838e-004 and the average number of iterations is 48038.
\label{BRTLambertianData1}}
\end{table}

These results, together with the results for specular reflection from \cite{L1}, show that the reconstruction error of the new numerical solution of 
the broken ray tomography problem using Lambertian reflection is approximately three times smaller compared to an average reconstruction error using 
specular reflection of 3.874848e-005 by the algorithm from \cite{L1} for the same test function, domain and obstacle and an order of magnitude smaller 
than the reconstruction error for classical tomography in the presence of a reflecting obstacle. On average, the new method also appears to be faster 
compared to classical tomography with ART. The order of the rays in both the ART and BRTL tests in this work is randomized before applying the Kaczmarz 
method therefore the speed of convergence difference is due to the use of reflection. We will report further results on the speed of convergence.

Table \ref{BRTLambertianData2} compares for the same environment and test function results on the accuracy of broken ray tomography and tomography 
without reflection when the size of the obstacle varies. Again, on average, broken ray tomography with Lambertian reflection is approximately 
an order of magnitude more accurate compared to classical tomography with ART. Table \ref{BRTLambertianData3} shows that for the same environment
and test function when the fraction of broken rays is decreased below a threshold the reconstruction error increases significantly. For example, 
for the data in Table \ref{BRTLambertianData3}, the error increases significantly when the number of broken rays is $10\%$ or less of the number of all rays. 
Below the threshold, the number of broken rays is so small that the tomography problem becomes approximately tomography without reflection. 
The set of boundary reflection points for the obstacle from Table \ref{BRTLambertianData3} excludes the four vertices of the square obstacle. 
The results from Table \ref{BRTLambertianData3} show that in practice tomography with reflection is robust with respect to the number of broken rays in the 
input set and remains in the same high-accuracy range, compared to the results from Table \ref{BRTLambertianData1} and Table \ref{BRTLambertianData2},
when some broken rays, e.g. no more than $90\%$  of all rays, or even when some boundary points are excluded from the computation or from the input data.
In the reflection model for Table \ref{BRTLambertianData3} and Table \ref{BRTLambertianData4} reflection from vertices is ignored and this result is 
applicable to reflection from edges in the solution of the tomography problem in $\mathbb{R}^3$.

Table \ref{BRTLambertianData4} presents data for BRTL and ART tomography with thirteen different test functions. The results from Table \ref{BRTLambertianData4} confirm previous results
that BRTL is approximately an order of magnitude more accurate than ART and has faster speed of convergence. 

\begin{table}[h]
\begin{tabular}{|c|cc|cc|}
	\hline
Side Length & ART Error  & ART Iterations & BRTL Error  & BRTL Iterations  \\
  \hline

130 & 1.753857e-004 & 48820 & 4.034727e-005 & 48341 \\
156 & 1.703602e-004 & 72874 & 2.385534e-005 & 51308 \\
182 & 1.665381e-004 & 94962 & 5.970411e-005 & 32424 \\
208 & 2.238438e-004 & 67407 & 3.354918e-005 & 38949  \\
234 & 2.599038e-004 & 53697 & 7.593578e-006 & 50307 \\
260 & 4.191898e-004 & 23842 & 2.235613e-005 & 36354 \\
286 & 1.663114e-004 & 113613 & 2.677853e-005 & 30901 \\
312 & 2.894623e-004 & 35652 & 5.080793e-005 & 24581 \\
338 & 1.292861e-004 & 97701 & 6.638226e-006 & 40774 \\
364 & 2.198555e-004 & 33167 & 8.823598e-006 & 36061 \\

  \hline

Average & 2.220137e-004 & 64173.5 & 2.804539e-005 & 39000.0 \\

  \hline
\end{tabular}
\caption{Error and number of iterations for broken ray tomography with Lambertian reflection(BRTL) at the boundary of the reflecting obstacle 
for a fixed number of 126050 rays with 50\% broken and 50\% unbroken rays.
Tomographic reconstruction is performed in ten experiments with different side lengths of the square obstacle.
The average error is 2.804539e-005 and the average number of iterations of the Kaczmarz method for finding a solution is
39000.0. The results for ART tomography with 126050 unbroken rays are shown in the left two columns of the table. 
The average error for ART is 2.220137e-004 and the average number of iterations is 64173.5.
\label{BRTLambertianData2}}
\end{table}

\begin{table}[h]
\begin{tabular}{|c|cc|}
	\hline
Fraction of Unbroken Rays & BRTL Error & BRTL Iterations  \\
  \hline

0.50 & 8.696576e-006 & 35457 \\ 
0.55 & 4.200643e-005 & 24092 \\
0.60 & 1.300616e-005 & 38860 \\
0.65 & 2.297882e-005 & 34489 \\
0.70 & 3.107772e-005 & 35637 \\
0.75 & 7.089215e-005 & 26228 \\
0.80 & 2.779947e-005 & 47781 \\
0.85 & 5.335306e-005 & 47160 \\
0.90 & 1.535794e-004 & 27335 \\
0.95 & 1.410281e-004 & 37328 \\

  \hline
\end{tabular}
\caption{Performance of broken ray tomography with Lambertian reflection for a fixed instance of the tomography problem with ray sets of 126050 rays with 
different fractions of broken and unbroken rays. When the fraction of unbroken rays is close to 1 the reconstruction error increases. The four vertices
on the square obstacle's boundary are excluded from the set of reflection points. 
\label{BRTLambertianData3}}
\end{table}

\begin{table}[h]
\begin{tabular}{|c|cc|cc|}
	\hline
Experiment & ART Error  & ART Iterations & BRTL Error  & BRTL Iterations  \\
  \hline

$f_0(x,y)$ & 1.887775e-004 & 84752 & 9.875051e-006 & 49777  \\ 
$f_1(x,y)$ & 2.108100e-004 & 72908 & 2.749901e-005 & 31266  \\ 
$f_2(x,y)$ & 2.022510e-004 & 79902 & 1.695669e-005 & 36142  \\ 
$f_3(x,y)$ & 1.187531e-004 & 39741 & 7.966167e-006 & 32711  \\ 
$f_4(x,y)$ & 1.604755e-004 & 102706 & 3.755106e-006 & 56232 \\ 
$f_5(x,y)$ & 1.525628e-004 & 115927 & 5.800111e-005 & 23606 \\ 
$f_6(x,y)$ & 2.737433e-004 & 51518 & 1.060598e-005 & 40013  \\ 
$f_7(x,y)$ & 3.425007e-004 & 36230 & 1.743366e-005 & 36216  \\ 
$f_8(x,y)$ & 1.711886e-004 & 142001 & 7.532997e-006 & 46364 \\ 
$f_9(x,y)$ & 2.788553e-004 & 61025 & 1.675307e-005 & 37604   \\ 
$f_{10}(x,y)$ & 1.973837e-004 & 101929 & 1.024641e-005 & 39846 \\ 
$f_{11}(x,y)$ & 2.764046e-004 & 49940 & 8.977682e-006 & 43783  \\ 
$f_{12}(x,y)$ & 2.725743e-004 & 45101 & 1.218634e-005 & 42030  \\ 
  \hline

Average &  2.189446e-004 & 75667.692308 & 1.598379e-005 & 39660.769231 \\

  \hline
\end{tabular}
\caption{Error and number of iterations for broken ray tomography with Lambertian reflection(BRTL) at the boundary of the reflecting obstacle 
for a fixed number of 126050 rays and thirteen different functions. The average error for BRTL is 1.598379e-005 and the average number 
of iterations for finding a solution is 39660.769231, and are shown in the right two columns of the table. 
The results for classical tomography with the ART method with the same number of rays are shown in the left two columns of the table. 
The average error for ART is 2.189446e-004 and the average number of iterations is 75667.692308.
\label{BRTLambertianData4}}
\end{table}

\section{Reduction of the Number of Equations of Overdetermined Linear Systems for the Tomography Problem}\label{drburzo}

We have shown that tomography with both specular and Lambertian reflection is more accurate than ART tomography in 
the presence of a reflecting obstacle. In this section, we present new algorithms and architectures for reducing
the number of equations of the overdetermined linear systems that are solved by the algorithms for broken ray tomography and 
other algebraic reconstruction methods. The reduced number of equations in linear systems for solving the tomography problem results in improved speed, 
i.e. smaller absolute running time, of the tomographic reconstruction algorithms from Section \ref{section_BRTL} and \cite{L1} 
as well as the speed of ART tomography and other algebraic reconstruction methods. 
 
In order to motivate the new approach, we return to specular reflection and consider the case when both the obstacle and observation 
boundary are reflective and the transmitters and receivers are transceivers. Then a ray sent from a transmitter $T_1$ can reflect at obstacle point $O_1$, reflect again at a receiver 
point $R_1$, and then depending on the geometry of the boundary reflect at a second obstacle point $O_2$ or directly reach receiver 
$R_2$ via an unbroken ray from $R_1$ to $R_2$, and so on. We consider only ray paths with multiple reflections that do not self-intersect and use these paths, 
instead of unbroken and broken rays composed of one or two segments, as the input ray paths for algebraic reconstruction algorithms. 

The basic computational model for algebraic reconstruction tomography methods is a computation grid composed of cells 
that are intersected by ray paths \cite{KS}. First, we extend this model to tomography in the presence of an obstacle: 
consider a square domain M in $\mathbb{R}^2$ that contains the observable domain $\Omega_0$. The square M is subdivided into a grid of $N^2$ 
squares or cells of size d. Define $f(x)=0$ in $\Omega_1$, inside the obstacle, and in $M \backslash \Omega_0$, outside the observation boundary, 
and look for a good approximation of $f(x)$ in each cell of the grid. The value of $f(x)$ is considered to be constant in each cell. 
We arrange linearly the $N^2$ cells into a column vector $f=(f_{1,1},..., f_{N,N})$ where $f_{i,j}$ is the value of $f(x)$ in cell $M[i][j]$.
When a ray j intersects cell i of the vector f, then the length of the ray segment that the cell cuts from the ray segment is the weight of the cell with respect to this ray 
or $w[j][i]$. When a cell is intersected more than once, we choose one of the weights for the cell or average their values in order to get one weight for the cell. The matrix of weights for all cells and all rays is denoted as W and has $N_r$ rows and $N^2$ columns, where $N_r$ is total number of rays. Let $P_j>0$ be the 
travel time of ray j. Let $$ \sum_{i=1}^{N^2} w[j][i]f[i] = P_j $$

Then the linear system of equations for all rays' travel-times can be expressed as $$Wf=P$$ where P is the column vector of ray travel times. 
In the algorithm from Section \ref{section_BRTL}, the matrix W is obtained from the set of rays or ray paths L and P is the vector of corresponding travel times for these rays.
The system $Wf=P$ is then solved by the Kaczmarz method:

$$ f^{(i+1)} = f^{(i)} + \frac{(P_{h}-w_{h} \cdot f^{(i)})}{w_{h} \cdot w_{h}}w_{h} $$

where $f^{(0)}$ is an initial guess, $h = (i \bmod N_r) + 1$ and $N_r$ is the number of rays or number of rows of the linear system.

In order to reduce the number of equations of $Wf=P$ we define abstract rays. 

\begin{definition}\label{ray}
An abstract ray in $\Omega_0$ is a finite set $r=\{\gamma_1, \gamma_2, ...,\gamma_n\}$ where $\gamma_1,\gamma_2, ..., \gamma_n$ are broken or unbroken rays
that do not intersect except at their endpoints. 
The travel time for r is $t = t_{1}+...+t_{n}$ where $t_i$ is the travel time for $\gamma_i$, $1 \leq i \leq n$.
\end{definition}

The above definition includes ray paths with multiple reflections that do not self-intersect and are allowed to have gaps.
Our algorithm for reducing the number of equations of $Wf=P$ is to partition a finite set $L$ of broken and unbroken rays into a set $L'$ of abstract rays.
Note that $ |L'| \leq |L| $ and when $L$  is sufficiently large e.g. $L=\mathbb{A}$, where $\mathbb{A}$ is a finite set of rays representing 
all discrete rays in $\Omega$ as in \cite{L1}, then there exists a partition of $L$ into a set $L'$ of abstract rays such that $$ |L'| < |L| $$ 
This result follows from the fact that a sufficiently large set of rays $L$ contains broken and unbroken rays that 
do not intersect and we can combine these rays into abstract rays which leads to a smaller number of abstract rays in $L'$ than the number of rays in $L$.  
Alternatively, instead of using this algorithm, if our scanner architecture 
delivers a set $L'$ of abstract rays e.g. by measuring multiple-reflections, then again $ |L'| < |L| $ for any set $L'$ that contains at least one abstract 
ray $r=\{\gamma_1, \gamma_2, ...,\gamma_n\}$ such that $n \geq 2$, i.e. $|r| \geq 2$. $L$ is the union of all broken and unbroken rays that
are contained in abstract rays r such that  $ r \in L'$. For such scanner architectures we could require, although it is not necessary, that an abstract ray is 
an ordered set and that its ray elements are ordered in time.

The set of abstract rays $L'$ and the set of corresponding travel times $P'$ for the abstract rays in $L'$ can be used directly in the algorithm
from Section \ref{section_BRTL}.  We construct $W'$ and $P'$ by the following procedure: for each abstract ray $r_j$
intersect the cells of the grid with all broken and unbroken ray elements $\gamma_u \in r_j$,$1 \leq u \leq n$, $|r|=n$ and assign weights $w[j][i]$,
where i is the number of the cell that is intersected in the unknown vector f. The weights are well defined except for
cells that are intersected more than once. For a cell that is intersected by several rays $\gamma_u$ define the weight of the cell as the average of the weights 
determined by the different ray segments that intersect the cell. Then each abstract ray in $L'$ determines one row from a weight matrix $W'$ with $|L'|$ 
rows and the corresponding element of $P'$ is the travel time t of the abstract ray. This time is given by Definition \ref{ray} as the sum of travel times of r's 
broken and unbroken ray elements.The error introduced by each cell that is intersected more than once is $$ e_c \leq \sqrt{2} \frac{d}{v_{min}} $$ where d is the length of the side of the square cell and 
$v_{min}>0$ is a lower bound on the velocity in the domain. When $d$ is sufficiently small this error term is very small.

The abstract ray representation results in fewer equations for $W'f=P'$ compared to $Wf=P$ because as discussed for sufficiently large L, $ |L'| < |L| $.
The presence of an obstacle and the requirements of the theory of broken ray tomography for considering all broken and unbroken rays increase 
significantly the size of the linear system that has to be solved. If the domain $\Omega_0$ is contained in a computation grid with size $N$, i.e. 
a square with $N^2$ square cells, we need a linear system with $O(N^2)$ equations and therefore $O(N^2)$ rays in order to solve  
the tomography problem via the Kaczmarz method without reflection. Position a transceiver, i.e. a transmitter and 
a receiver, in the center of each cell side of the grid's boundary. The total number of transceivers $U$ is then 4N. Consider a discretization of the obstacle's boundary that provides 
$B$ reflection points. Then we have that for specular reflection the number of all broken and unbroken rays $n_{au}+n_{ab}$ is $O(U^2+UB)$. For example, 
for a grid with 4096 cells and tomography without reflection we need a linear system of approximately 4096 equations provided by 4096 measurement rays. 
For tomography in the presence of obstacles using the same grid and set of transmitters and receivers, we approximate the infinite set of all broken and ubroken rays 
with 126050 rays which leads to an overdetermined linear system with 126050 equations and 4096 unknowns. In order to speed up the solution of the system we can transform it into
a system with a much smaller number of equations where each equation of the new system is determined by an abstract ray with zero, one or more reflection points. 
The equations from the original system determined by the input set of broken and unbroken rays are 
replaced by a system with fewer equations determined by fewer abstract rays containing the input broken and unbroken rays. This approach incorporates information from all broken and unbroken rays into the
linear system and reduces the size of the linear system. These pre-processing steps can be parallelized and encapsulated in the following algorithm for constructing $L'$:

S1. Start with an input set $L=\mathbb{I}$ that contains a discrete approximation of the set of all broken and unbroken rays in $\Omega_0$.

S2. Pick any broken or unbroken ray $r_1 \in \mathbb{I}$. Let $L_1=\{r_1\}$. 

S3. Pick any ray $r_2$ from the remaining rays in $\mathbb{I}$ such that $r_1$ and $r_2$ have a common vertice on the observation boundary 
$\partial{\Omega_0}$. This implies that the transmitter of $r_2$ is a receiver of $r_1$ or that the transmitter of $r_1$ is a receiver of $r_2$. 
For specular reflection and for generating abstract ray paths without gaps we make the additional check that at vertice P shared by both rays, the ray segments from $r_1$ and $r_2$ that intersect at P
form angles of $\alpha$ and $\pi-\alpha$ with the tangent at P. In other words, for specular reflection and ray paths without gaps, we check that for the concatenated ray
reflection at P is mirror-like. Let $L_1 = \{r_1, r_2\}$.

S4. Repeat step S3 as many times as possible by adding to the beginning or end of $L_1$ a new ray $r$ from the remaining rays in $\mathbb{I}$ such
that $L_1$ does not self-intersect. Stop when $L_1$ can not be extended by a new ray $r$ and remain non self-intersecting.   

S5. Repeat steps S2-S4 for the remaining rays in $\mathbb{I}$ and obtain rays $L_2$, $L_3$,...,$L_n$ such that each ray $L_i$ is defined as in Definition \ref{ray}.   

The travel time for ray $L_i$ is the sum of travel times of the individual rays from which $L_i=\{r_{i1}, ..., r_{is}\}$ is composed:
$$t_i = t_{i1}+...+t_{is} $$

S6. $L'= \{L_1,...,L_n\}$ contains all abstract rays constructed in the previous steps.

The Kaczmarz method allows rays with flexible geometries that do not self-intersect. We construct $W'f=P'$ and
solve it by the Kaczmarz method. This provides a faster version of the algorithm from section \ref{section_BRTL} 
with input of $L'=\{L_1$, $L_2$,...,$L_n\}$. The method is applicable to other algebraic reconstruction algorithms including ART
because their input ray sets can be partitioned into abstract rays.

The reduction of the number of equations of the linear system for solving the tomography problem depends on $n$ which in turn depends on the
number of rays that compose each individual abstract ray. For example, if the average number of rays that compose each ray in $L'$ is 2 then then number
of equations in the linear system will be reduced by half. An alternative to this method for partitioning the input set of broken and unbroken
rays is to scan with rays with multiple reflections. Consider a scanner architecture with a reflective observation boundary as described above.
When a transceiver sends a ray $L_i$ it could reflect several times and is an abstract ray. The scanner keeps a bit map that corresponds to the
approximated set $\mathbb{A}$ of all broken and unbroken rays in $\Omega_0$. If a ray r from the set $\mathbb{A}$ is a subset of $L_i$ then the 
corresponding bit for r is set to true. The scanning process stops if all bits in the bit map are set to true. 
Using integer coordinates and a data structure such as a hash table, the lookup operation for a broken or unbroken ray from $L_i$ is $O(1)$. 
Therefore, the performance of algebraic reconstruction with abstract rays is determined by the time for solving the linear system $W'f=P'$
and is faster compared to methods for solving the linear system $Wf=P$.

The method for reducing the number of equations in the linear system for broken ray tomography provides a flexible framework for fast 
tomography in the presence of obstacles. It is also applicable to improving the performance of a wide class of algebraic 
reconstruction methods.

\section{Acknowledgements}

I would like to thank Professor Gregory Eskin for his continuous guidance. 

\section*{References}

\bibliography{BRT}

\end{document}